\input amstex.tex
\documentstyle{amsppt}
\NoRunningHeads 
\magnification\magstep1
\def\mpr#1{\;\smash{\mathop{\hbox to 20pt{\rightarrowfill}}\limits^{#1}}\;} 
\def\mpl#1{\;\smash{\mathop{\hbox to 35pt{\leftarrowfill}}\limits^{#1}}\;} 
\def\mpd#1{\big\downarrow\rlap{$\vcenter{\hbox{$\scriptstyle#1$}}$}}
\def\mdp#1{\llap{$\vcenter{\hbox{$\scriptstyle#1$}}$}\big\downarrow}
\def\mpu#1{\big\uparrow\rlap{$\vcenter{\hbox{$\scriptstyle#1$}}$}}
\def\Qt{\Bbb Q}

\def\twoheadrightarrow{\longrightarrow\hskip -20pt\longrightarrow}
\def\X{{\widetilde{X}}}
\def\Y{{\widetilde{Y}}} 

{\nopagenumbers
\topmatter
\title A MORPHISM OF INTERSECTION HOMOLOGY
\\ INDUCED BY AN ALGEBRAIC MAP \endtitle
\author Andrzej Weber \endauthor

\thanks I would like to acknowledge the hospitality of
the Institut de Math\`ematiques
de Luminy and to thank Prof. J. P. Brasselet
personally for he encouraged me to find a simpler proof of the main theorem 
from \cite{BBFGK}.
I would also like to thank the Authors of \cite{BBFGK} for correcting 
this note. 
\endthanks
\thanks Partially supported by KBN 2P30101007 grant\endthanks
\affil Institute of Mathematics, University of Warsaw \endaffil
\address ul. Banacha 2, 02--097 Warszawa, Poland \endaddress
\email aweber\@mimuw.edu.pl \endemail

\abstract Let $f:X @>>> Y$ be a map of algebraic varieties. 
Barthel, Brasselet, Fieseler, Gabber and Kaup have shown that there
exists a homomorphism of intersection homology groups
$f^*:IH^*(Y)@>>>IH^*(X)$ compatible with the induced homomorphism on
cohomology; \cite{BBFGK}. The crucial point in the argument is
reduction to the finite characteristic. We give an alternative and
short proof of the existence of a homomorphism $f^*$. Our
construction is an easy application of the Decomposition Theorem.
\endabstract

\endtopmatter}
\vfill\eject
\pageno 1
\document

Let $X$ be an algebraic variety, $IH^*(X)=H^*(X;IC_X)$ its rational 
intersection homology group with respect to the middle perversity 
and $IC_X$
the intersection homology sheaf which is  an  object  of  derived 
category of 
sheaves over $X$; \cite{GM1}. We have the homomorphism $\omega_X:H^*(X;\Qt)\mpr{} 
IH^*(X)$ 
induced by the canonical morphism of the sheaves $\omega_X:\Qt_X\mpr{}IC_X$.

Let $f:X \longrightarrow Y$ be a map of
algebraic varieties. It induces a homomorphism of the cohomology groups. The 
natural question arises: Does there exist an induced homomorphism for
intersection homology compatible with $f^*\;$?
\vskip 1pt
$$\CD IH^*(Y) & \mpr{?} & IH^*(X)\\
\mpu{\omega_Y} & & \mpu{\omega_X} \\
H^*(Y;\Qt) & \mpr{f^*} & H^*(X;\Qt)\,.
 \endCD$$
This question has a positive answer. 
The authors of \cite{BBFGK} proved the following:
 
\proclaim
{Theorem 1} Let $f:X \longrightarrow Y$ be an algebraic map of
algebraic varieties. Then there exists a morphism $\lambda_f: IC_Y 
\longrightarrow
Rf_*IC_X$ such, that the following diagram with the canonical morphisms 
commutes: 
\vskip 1pt
$$\CD IC_Y & \mpr{\lambda_f} & Rf_*IC_X\\
\mpu{\omega_Y} & & \mpu{Rf_*(\omega_X)} \\
\Qt_Y & \mpr{\alpha_f} & Rf_*\Qt_X\,.
\endCD$$
\endproclaim 

In fact, \cite{BBFGK} proves the existence of a morphism
$\mu_f: f^*IC_Y\rightarrow IC_X$, which is adjoined to $\lambda_f$.

The sheaf language can be translated to the following: 
an induced homomorphism of intersection homology exists in a functorial way with 
respect to the open subsets of $Y$. This means that there exists a 
compatible family of induced homomorphisms
$$f^*_{\lambda,U}:IH^*(f^{-1}U)\mpr{}IH^*(U)\,,$$
which is also compatible with the family
$$f^*_{|U} :H^*(f^{-1}U;\Qt)\mpr{}H^*(U;\Qt)\,.$$
As shown in \cite{BBFGK} the morphism $\lambda_f$ (and $\mu_f$)
is not unique. It is not possible to choose the morphisms
$\lambda_f$ (nor $\mu_f$) in a functorial way with respect to all algebraic 
maps (p.160). The simplest counterexample is the inclusion
$\{(0,0)\}\hookrightarrow \{(x_1,x_2):\, x_1x_2=0\}$, which
can be factored through the inclusions 
$\{(0,0)\}\hookrightarrow \{(x_1,x_2):\, x_i=0\}$ for $i=1$ or 2.

We will give a short proof of the main theorem from \cite{BBFGK}.
We will derive it from the Decomposition Theorem. The
refe\-rence to the Decomposition Theorem is \cite{BBD, 6.2.8}
(see also \cite{GM2}) and slightly more general \cite{Sa}.

We will use only the following corollary from the Decomposition
Theorem:

\proclaim
{Corollary from the Decomposition Theorem} Let $\pi:X
\longrightarrow Y$ be a proper surjective
map of algebraic varieties. Then $IC_Y$ is a direct summand in
$R\pi_*IC_X$.
\endproclaim

The idea of the proof of our Theorem is simple, the essence is the
argument similar to \cite{BBFGK, Remarque pp.172--174}. We take a
resolution $\pi_Y:\Y\longrightarrow Y$ and enlarge the space $X$ to
obtain a map $\tilde f:\X \longrightarrow \Y$. 
There exists the induced morphism of intersection
homology $\lambda_{\tilde f}$ for $\tilde f$. By the Decomposition 
Theorem the intersection homology of $X$ (and $Y$) is a direct
summand of intersection homology of $\X$ (resp. $\Y$). We compose
$\lambda_{\tilde f}$ with the projection and inclusion in the direct
sums to obtain the desired morphism $\lambda_f$.
\vskip 3pt

\remark{Remark} If we insisted then $\X$ might be even smooth
of the same dimension as $X$ with the map $\pi_X:\X \longrightarrow X$
generically finite; compare \cite{BBFGK, p.173}. \endremark

\demo{Proof of Theorem 1} We may assume that $X$ and $Y$ are irreducible.
Let $\pi_Y:\Y\longrightarrow Y$ be a resolution of Y. Denote by $\X$
the fiber product (pull-back) $X\times_Y\Y$. Note that it is a
variety, which may be singular and not equidimensional. We
have a
commutative diagram of algebraic maps ($\pi_X$ and $\pi_Y$ proper):
$$\CD \X & \mpr{\tilde{f}} & \Y \\
\mdp{\pi_X} & & \mpd{\pi_Y} \\
X & \mpr{f} & Y \endCD$$ 
and the associated diagram of sheaves over $Y$:
\vskip 1pt
$$\CD R\pi_{Y*}IC_\Y & = & R\pi_{Y*}\Qt_\Y & 
\;\mpr{R\pi_{Y*}(\alpha_{\tilde f})}\; & Rf_*R\pi_{X*}\Qt_\X &
\; \mpr{Rf_*R\pi_{X*}(\omega_\X) } \; & Rf_*R\pi_{X *}IC_\X 
\\
?\big\uparrow & & \mpu{\alpha_{\pi_Y}} & & \mpu{Rf_*(\alpha_{\pi_X})} & & 
\big\downarrow ? \\
IC_Y &\mpl{\omega_Y }& 
\Qt_Y & \mpr{\alpha_f} & Rf_*\Qt_X & \mpr{\quad Rf_*(\omega_X)\quad }&
Rf_*IC_X \,. 
\endCD
$$ 
To prove the existence of a morphism $\lambda_f:IC_Y\longrightarrow Rf_*IC_X$, 
we will show that the arrows with question marks exist
in a way that the diagram remains commutative. The existence of such
morphisms follows from the Decomposition Theorem for $\pi_Y$ and
$\pi_X$, see Corollary. 
The sheaf $IC_Y$ is a direct summand in $R\pi_{Y*}IC_{\Y}$:
$$i:IC_Y\hookrightarrow R\pi_{Y*}IC_{\Y}\,.$$
We also have a projection:
$$p:R\pi_{X*}IC_\X \twoheadrightarrow IC_X\,,$$
which induces
$$Rf_*(p):Rf_*R\pi_{X*}IC_\X \twoheadrightarrow Rf_*IC_X\,.$$
It remains to prove the commutativity of the diagram. We compare the
morphisms over $Y$:
$$ \Qt_Y \mpr{\omega_Y } IC_Y \mpr{i} R\pi_{Y*}IC_{\Y}=R\pi_{Y*}\Qt_\Y$$
and the natural one
$$ \Qt_Y \mpr{\alpha_{\pi_Y}} R\pi_{Y*}\Qt_\Y\,.$$
Respectively over $X$ we compare the morphisms:
\vskip 1pt
$$ \Qt_X \mpr{\alpha_{\pi_X}} R\pi_{X*}\Qt_\X
\mpr{\;\;R\pi_{X*}(\omega_\X)\;\; } R\pi_{X*}IC_\X 
\mpr{p}
IC_X$$ and the canonical one
$$ \Qt_X \mpr{\omega_X} IC_X\,.$$
Let $U$ (resp. $V$) be the regular part of $Y$ (resp. $X$).
After multiplication by a constant if necessary, these morphisms are
equal on $U$ (resp. on $V$).
We will show that an equality of morphisms over an open set implies the
equality over the whole space. We have the restriction morphism
$$\CD Hom(\Qt_Y, R\pi_{Y*}\Qt_\Y) & \mpr{\rho_U} & Hom((\Qt_Y)_{|U},
(R\pi_{Y*}\Qt_\Y)_{|U}) \\
\| & & \| \\
H^0(\Y) & & H^0(\pi_Y^{-1}(U))\,.\endCD$$
The kernel of $\rho_U$ is
$H^0(\Y,\pi^{-1}_Y(U))$, which is trivial. 
We have the same for the morphisms over $X$:
$$\CD Hom(\Qt_X, IC_X) & \mpr{\rho_V} & Hom((\Qt_X)_{|V},(IC_X)_{|V}) \\
\| & & \| \\
IH^0(X) & & IH^0(V)\,.\endCD$$
The kernel is $IH^0(X,V)=0$.\qed
\enddemo

\remark{Remark} The restriction morphisms $\rho_U$ and $\rho_V$ are in fact 
isomorphisms. The cokernel of $\rho_U$ is contained in
$H^1(\Y,\pi^{-1}_Y(U))=H^{cl}_{2(\dim \Y)-1}(\Y\setminus\pi^{-1}_Y(U))$
which is trivial for dimensional reason. The second follows from
\cite{Bo, V.9.2 p.144} as noticed in \cite{BBFGK, p.178}.\endremark

\enddocument
\end